\documentclass[draft,leqno]{amsproc}
\usepackage{latexsym,amssymb,amsmath,color,pifont}

\newtheorem{thm}{Theorem}
\newtheorem{lem}[thm]{Lemma}

\newtheorem{pro}[thm]{Proposition}
\newtheorem{cor}[thm]{Corollary}

\theoremstyle{remark}
\newtheorem{rem}[thm]{Remark}

\DeclareMathOperator{\lin}{lin}

\newcommand*{\cbb}{\mathbb C}

\newcommand*{\D}{\mathrm d}

\newcommand*{\I}{\mathrm i}

\newcommand*{\is}[2]{\langle#1,#2\rangle}

\newcommand{\mfrak}{\mathfrak m}
\newcommand*{\nbb}{\mathbb N}

\newcommand*{\nul}{\{0\}}
\newcommand*{\okr}{\stackrel{\scriptscriptstyle{\mathsf{def}}}{=}}
\newcommand*{\ol}{\overline}
\newcommand*{\pcal}{\mathcal P}

\newcommand*{\rbb}{\mathbb R}

\newcommand*{\tbb}{\mathbb T}
\newcommand{\trans}{{\mathop\intercal}}
\newcommand*{\zbb}{\mathbb Z}

\newcommand{\emv}{\stackrel{{\scriptscriptstyle{{\mathit V}}}}{=}}
\newcommand{\emvdwa}{\stackrel{{\scriptscriptstyle{{\mathit
V_2}}}}{=}}

\begin{document}
   \title[Christoffel-Darboux kernels in several real variables]{Christoffel-Darboux kernels in several real variables}
   \author[D. Cicho\'n, F.H. Szafraniec]
{Dariusz Cicho\'n, Franciszek Hugon
Szafraniec}
   \address{Wydzia\l\ Matematyki i Informatyki, Uniwersytet
Jagiello\'nski, ul.\ {\L}ojasiewicza 6, PL-30348
Kra\-k\'ow}
   \email{dariusz.cichon@uj.edu.pl}
   \email{umszafra@cyfronet.pl}
   \begin{abstract}
The Christoffel-Darboux kernels for orthogonal polynomials in several real variables are investigated within the context of the three term recurrence relation reformulated for this purpose. Examples of orthogonal polynomials on the unit circle and on the Bernoulli lemniscate are discussed.
   \end{abstract}
   \maketitle
The recent progress in orthogonal polynomials in several variables entails the three term recurrence relation, which in this case requires organizing the polynomials in columns and employing the matrix coefficients in place of ordinary numbers. Results of this type were obtained in \cite{kow1} and \cite{xu1}, and later in \cite{css} where the authors managed to include the case of orthogonal polynomials with respect to measures with thin support, where by ``thin'' we mean ``contained in a proper algebraic subset of $\rbb^n$.'' Knowing the three term recurrence relation we may formulate criteria for existence of a Borel measure, orthonormalizing the system of polynomials. Another kind of application is possibility of writing the Christoffel-Darboux formula, which is well-known for polynomials in a single real variable, and highly expected to hold in several variables too. Our paper answers this call in the way it seems to be possible in the context of polynomials in several variables---and not only possible though intrinsically reasonable. Since the equations considered here hold mostly modulo an ideal it is impossible to perform division by a polynomial, as this happens in the case for one variable polynomials. This justifies our approach with no division involved. One of surprising effects is that the Christoffel-Darboux formula requires equations modulo another ideal in twice as many variables, which means ``doubling'' the starting ideal. This has made us consider some abstract ideas which turn out to be indisposable in the process of deriving the Christoffel-Darboux formula in question. The usual way is to go from the three term recurrence relation to the Christoffel-Darboux formula, and this has motivated us to consider some examples leading to the relation in question (or, more specifically, showing the possible way of obtaining it recursively). Among the examples there are polynomials which are orthogonal with respect to a measure supported in the unit circle or the Bernoulli lemniscate.

Our paper may be considered a challenge to the opinion appearing in \cite[p.\ 299]{simon1} that ``the conventional wisdom is that there is no CD formula for general OPs''. This however is less important motivation, as our primary one is to draw attention to the Christoffel-Darboux formula itself, which surprisingly reduces an expression involving several initial orthogonal polynomials to the formula with only two last ones (except that for several variables things get more complicated).

A nice by-product we include is a very concise version of a theorem stated in \cite{css} which here appears as Theorem \ref{gfav}. The references we list at the end could be easily made at least twice as long so if needed, the interested reader is encouraged to confront the references in \cite{css}.

\section{Introduction}
Let $\nbb=\{0,1,2,\ldots\}$. We make use of the Kronecker delta $\delta_{ab}$ which is equal to $1$ if $a=b$ or $0$ otherwise, where $a$ and $b$ are any objects. All measures on $\rbb^d$ are assumed to be Borel measures with all moments finite, i.e.\ $\int_{\rbb^d} \|x\|^n \D\mu(x) <\infty$ for all $n\geqslant 0$. Let $\pcal_d$ stand for the space of all complex polynomials in $d$ (real) variables, and $\pcal_d^k$ -- for the space of all complex polynomials in $d$ variables of degree at most $k$. We say that a measure $\mu$ on $\rbb^d$ {\em orthonormalizes} a sequence $\{p_k\}_{k=0}^\infty \subset \pcal_d$, if $\int_{\rbb^d} p_k \bar p_l \D\mu = \delta_{kl}$ for all $k,l\in \nbb$, where $\bar p$ is the defined by $\bar p(x)=\overline{p(x)}$, $x\in\rbb^d$, $p\in\pcal_d$.

Let us recall the Favard theorem on orthogonal polynomials in one variable. The following version of this theorem is a direct consequence of Theorem 53 in \cite{css} (cf.\ \cite{chih}).
\begin{thm}\label{fav}
   If $\{p_k\}_{k=0}^\infty$ is a sequence of real polynomials in one variable such that $p_0 = 1$ and $\deg p_k=k$ for all $k\in\nbb$, then the following two conditions are
   equivalent:
   \begin{enumerate}
   \item[(i)] there exists a
measure $\mu$ on $\rbb$ which orthonormalizes
$\{p_k\}_{k=0}^\infty$,
   \item[(ii)] for every $k \in \nbb$, there exist $a_k \in
\rbb$ and $b_k \in \rbb$ such that
   $$X p_k = a_k p_{k+1} + b_k p_k +  a_{k-1}
p_{k-1}, \text{ where } a_{-1} \okr 1 \text{ and } p_{-1} \okr 0.
   $$
   \end{enumerate}
   \end{thm}
The condition (ii) is customarily called {\em the three term recurrence relation}.

When considering a multi-variable version of this theorem the succesful attempts were undertaken by Kowalski \cite{kow1} and then further developed by Xu \cite{xu1}, yet they formulated their versions only for full polynomial bases of $\pcal_d$. This clearly excludes some interesting measures, e.g.\ the Lebesgue measure on the unit circle in $\rbb^2$, as in this case every polynomial is orthogonal to $x^2+y^2-1$. It was discovered in \cite{css} that this deficiency can be overcome by introducing equality modulo an ideal instead of dealing with ordinary equality in the three term recurrence relation. The aforesaid ideal in the case of the circle would be the set of all polynomials vanishing on the circle.

We now proceed to introduce our \underbar{main} notions required to state the Favard theorem in several variables. Let $V \subset \pcal_d$ be a proper ideal and $\Pi_V:\pcal_d \to \pcal_d/V$ be the canonical embedding. We additionally assume that $V$ is a $*$-ideal, i.e.\ $\overline p \in V$ whenever $p\in V$.
We will say that a set of polynomials is {\em linearly independent over $V$} (or {\em $V$-linearly independent}) if $\Pi_V$ maps this set injectively onto a linearly independent subset of $\pcal_d/V$. 
Let 
$$
d_V(k) \okr \dim \Pi_V(\pcal^k_d) - \dim \Pi_V(\pcal^{k-1}_d), \quad k\geqslant 1,
$$
and $d_V(0)\okr 1$. In most cases $d_V(k) \geqslant 1$ for all $k\geqslant 0$, in which case we set $\varkappa_V = \infty$, and otherwise $\varkappa_V = \max\{k\geqslant 1 \colon d_V(k) \neq 0\}$; this definition makes sense because by formula (8) in \cite{css} $d_V(k+1) = 0$ provided $d_V(k)=0$, $k\geqslant 1$.  Furthermore, a sequence $\{Q_k\}_{k=0}^{\varkappa_V}$ is called a {\em rigid $V$-basis} of $\pcal_d$ if every $Q_k$ is a column polynomial of size $d_V(k)$, i.e.
\begin{align}\label{kuka}
Q_k = \begin{bmatrix}q_1^{(k)} \\ \vdots \\ q_{d_V(k)}^{(k)}\end{bmatrix},\quad q^{(k)}_j \in \pcal_d\ (j=1,\ldots, d_V(k)),
\end{align}
all polynomials in $Q_k$ are of degree $k$ and the set
$$
\{q^{(k)}_{j_k}+V\colon k\in\nbb,\  k\leqslant \varkappa_V,\ j_k=1,\ldots, d_V(k)\}
$$
is a basis of $\pcal_d/V$. As shown in \cite{css} such bases always exist, and moreover, they can be build from the monomials. If $p,q \in \pcal_d$, then the notation $p\emv q$ means that $p+V=q+V$ (or equivalently $p-q \in V$). If $P$ and $Q$ are column polynomials, then we write $P\emv Q$ if the columns are of the same size and all the entries of $P-Q$ are in $V$. If $\{Q_k\}_{k=0}^{\varkappa_V}$ is a rigid V-basis of $\pcal_d$, $p\in \pcal_d$ and $\deg p=k \leqslant \varkappa_V$, then there are unique scalar row vectors $C_0,\ldots,C_k$ such that $p \emv \sum_{j=0}^k C_j Q_j$ (this remains true if $k>\varkappa_V$ provided the summation is done over $j=0,\ldots,\varkappa_V$).

Given a linear functional $L:\pcal_d\to \cbb$ we write
$$
L([p_{k,l}]_{k=0}^m{}_{l=0}^n) = [L(p_{k,l})]_{k=0}^m{}_{l=0}^n,
$$
where $p_{k,l} \in \pcal_d$. This way we can make sense of $L(PQ^*)$, where $P$ and $Q$ are column polynomials (not necessarily of the same size), and
\begin{align*}
Q^* = \begin{bmatrix}\bar q_1 & \ldots & \bar q_m\end{bmatrix} \quad \text{if} \quad Q= \begin{bmatrix} q_1 \\ \vdots \\ q_m \end{bmatrix}.
\end{align*}
Finally, we say that $L$ {\em orthonormalizes} a rigid $V$-basis $\{Q_k\}_{k=0}^{\varkappa_V}$ of real polynomials\footnote{\;i.e.\ each column polynomial is composed of polynomials with real coefficients}, if $L(Q_kQ_l^\trans) = 0$ for all $k,l \in \nbb$, $k,l\leqslant \varkappa_V $, $k\neq l$, and $L(Q_kQ_k^\trans) = I$, $k\in\nbb$, $k\leqslant \varkappa_V$, where $I$ stands for the identity matrix of appropriate size (which in these settings is equal to $d_V(k)$). A linear functional $L:\pcal_d \to \cbb$ is called {\em positive definite} if $L(p\overline{p})\geqslant 0$ for all $p\in\pcal_d$.

We are now in a position to state our generalization of the Favard theorem which is in flavour of \cite{css}.
\begin{thm}\label{gfav}
Let $V\subset\pcal_d$ be a proper $*$-ideal and $\{Q_k\}_{k=0}^{\varkappa_V}$ be a rigid $V$-basis of real polynomials with $Q_0=1$. Then the following conditions are equivalent:
\begin{enumerate}
\item[(A)] there exists positive definite $L:\pcal_d\to\cbb$ which orthonormalizes $\{Q_k\}_{k=0}^\infty$ and such that $V\subset \ker L$;
\item[(B)] for every $j=1,\ldots,d$ there exist systems of scalar matrices $\{A_{k,j}\}_{k=0}^\infty$ and $\{B_{k,j}\}_{k=0}^\infty$ of appropriate size such that
\begin{align}\label{3term}
X_jQ_k \emv A_{k,j}Q_{k+1}+ B_{k,j}Q_k + A_{k-1,j}^* Q_{k-1}, 
\end{align}
for all $k\in\nbb$, $k\leqslant \varkappa_V$, and $j=1,\ldots,d$, where $A_{-1,j}=1$ and $Q_{-1}=0$; if $\varkappa_V < \infty$, then $A_{\varkappa_V,j}$ is a column $[1, \ldots, 1]^\trans$ with $d_V(\varkappa_V)$ entries and $Q_{\varkappa_V+1} = 0$.
\end{enumerate}
\end{thm}
\begin{proof} This is a specified version of Theorem 36 in \cite{css}, if conditions (A) and (B) therein are simplified by our assumptions that  $\{Q_k\}_{k=0}^{\varkappa_V}$ is a rigid $V$-basis of $\pcal_d$; in particular the lengthy condition (B) therein reduces to (B-i) only. 
\end{proof}
\begin{rem}\label{inj}
As follows from Theorem 36 in \cite{css}, the matrix $\begin{bmatrix}A_{k,1}^* & \ldots & A_{k,d}^* \end{bmatrix}^*$ must necessarily be injective for all $k\geqslant 0$. Even more, all the matrices $A_{k,j}$ and $B_{k,j}$ must be real. To see this one can take adjoints\footnote{\ the adjoint of a complex matrix $[a_{kj}]_{k,j=1}^n$ is given by $[\overline{a_{jk}}]_{k,j=1}^n$} of both sides in \eqref{3term} and then transpose them to deduce that \eqref{3term} holds with $(A_{k,j}^*)^\trans$ and $(B_{k,j}^*)^\trans$ in place of $A_{k,j}$ and $B_{k,j}$, respectively, which employs the assumption that all $Q_k$ are real column polynomials. Since $\{Q_k\}_{k0}^{\varkappa_V}$ is a rigid $V$-basis, the matrices in \eqref{3term} are unique, thus $(A_{k,j}^*)^\trans = A_{k,j}$ and $(B_{k,j}^*)^\trans = B_{k,j}$, which proves our claim. In particular, \eqref{3term} can be written with $A_{k-1,j}^\trans$ in place of $A_{k-1,j}^*$. Actually, all the matrices $B_{k,j}$ are symmetric. Indeed, if we multiply both sides of \eqref{3term} by $Q_k^\trans$ from the left, and then apply $L$ appearing in Theorem \ref{gfav} to the equation, after omitting the zero terms and noticing that $L(B_{k,j}Q_kQ_k^\trans) = B_{k,j} L(Q_kQ_k^\trans)$, we get
$$
L(X_jQ_kQ_k^\trans) = B_{k,j}, \quad k\in\nbb, j \in \{0,\ldots,d\}.
$$
It now suffices to justify the symmetry of the left-hand side:
$$
L(X_jQ_kQ_k^\trans)^\trans = L \big( (X_jQ_kQ_k^\trans)^\trans \big) = L \big( X_j (Q_kQ_k^\trans)^\trans \big) = L( X_jQ_kQ_k^\trans).
$$
Note that this idea provides another way of showing that the matrices $A_{k,j}$ and $B_{k,j}$ are real.
\end{rem}
As is easily seen the case when $\varkappa_V<\infty$ (in other terms $\dim \pcal_d/V < \infty$) corresponds to a finite number of linearly independent polynomials over $V$, thus in this case we can have only a finite number of orthogonal polynomials whenever $V\subset \ker L$. Indeed, one can readily verify that under the latter assumption any set $C\subset \pcal_d$ such that $L(p\bar q) = \delta_{pq}$, $p,q\in C$, must be $V$-linearly independent, thus its cardinality is no greater than $\dim \pcal_d/V$. 

In the single variable case Theorem \ref{gfav} leads to the Favard theorem (Theorem \ref{fav}) due to the well-known fact that every positive definite linear functional on $\pcal_1$ can be represented as an integral with respect to a nonnegative Borel measure on $\rbb$ (see the discussion concerning the statement (50) in \cite{css}, see also \cite{shoh} and \cite{chih}). Disappointingly, albeit challengingly so, the positive definiteness does not guarantee the existence of such a representing measure in the several variable case (again see the comment on (50) in \cite{css}). A substantial part of \cite{css} is devoted to conditions ensuring existence of representing measures in the multi-variable case.
\section{Christoffel-Darboux kernel in the case of real variables}\label{CDreal}
Let $V\subset\pcal_d$ be a proper $*$-ideal and let $\{Q_k\}_{k=0}^{\varkappa_V}$ be a rigid $V$-basis of $\pcal_d$. For every $n\in\nbb$, $n\leqslant \varkappa_V$, we define the associated Christoffel-Darboux kernel by the formula
\begin{equation}\label{cdk}
K_n(x,y)\okr\sum_{k=0}^n Q_k^*(y)Q_k(x), \quad x,y\in \rbb^d. 
\end{equation}
We use the shorthand notation $\bigcup_{k=0}^n Q_k$ for the set of all entries of all column polynomials $Q_k$, $k=0,\ldots,n$. It is worth noticing that each $K_n$ enjoys the reproducing property in the following sense (cf. formula (1.2.37) in \cite{simonbook1}).
\begin{pro}\label{repker}
If $\{Q_k\}_{k=0}^{\varkappa_V}$ is a rigid $V$-basis of $\pcal_d$ which is orthonormal with respect to an inner product $\is \cdot{\text{-}}$ in $\pcal_d$, and $X_n \okr \mathrm{lin}\bigcup_{k=0}^n Q_k$ for $n\in \nbb$, $n\leqslant \varkappa_V$, then $K_n$ is a reproducing kernel for $X_n$, which means that for every $p\in X_n$ we have the equation
$$
p(y) = \is p{K_n(\cdot,y)},\quad y\in\rbb^d.
$$
\end{pro}
\proof
Let us write $\bigcup_{k=0}^n Q_k = \{q_0,\ldots,q_m\}$. The polynomial $p\in X_n$ may be written uniquely as $p=\sum_{k=0}^m a_k q_k$ with some $a_k \in\cbb$. Similarly, $K_n(x,y) = \sum_{k=0}^m \ol{q_k(y)} q_k(x)$, $x,y\in\rbb^d$. Then
\begin{align*}
\is p{K_n(\cdot,y)} = \Big\langle \sum_{k=0}^m a_k q_k, \sum_{k=0}^m \ol{q_k(y)} q_k\Big\rangle = \sum_{k=0}^m a_k q_k(y) = p(y),
\end{align*}
which is the desired conclusion.
\endproof
If  $H_1,H_2$ are linear subspaces of $\pcal_d$, then we identify the algebraic tensor product $H_1\otimes H_2$ with the linear space spanned by all polynomials of the form $h_1\otimes h_2 \in \pcal_{2d}$, where
$$
(h_1\otimes h_2)(x,y) = h_1(x)h_2(y),\quad x,y\in \rbb^d.
$$

We may now formulate a multi-variable version of the Christoffel-Darboux formula.
\begin{thm}\label{CDformula}
Let $V\subset\pcal_d$ be a proper $*$-ideal, and $\{Q_k\}_{k=0}^{\varkappa_V}$ be a rigid $V$-basis of real polynomials with $Q_0=1$. Assume that $\{Q_k\}_{k=0}^{\varkappa_V}$ obeys the three term recurrence relation {\em (B)} in Theorem \ref{gfav}. Then $(x_j-y_j)K_n(x,y)$ is equal to
\begin{align}\label{Ktilde}
[A_{n,j}Q_{n+1}(x)]^\trans Q_n(y) - Q_n(x)^\trans [A_{n,j}Q_{n+1}(y)]
\end{align}
modulo the ideal $V_2 \okr V\otimes\pcal_d+\pcal_d\otimes V \subset \pcal_{2d}$ for all $n\in\nbb$, $n\leqslant \varkappa_V$, and $j=1,\ldots,d$.
\end{thm}
\begin{proof}
Fix $j\in \{1,\ldots, d\}$ and $n\in\nbb$, $n\leqslant \varkappa_V$. Let the expression \eqref{Ktilde} be denoted by $\tilde K_n(x,y)$. Applying (B) supported by Remark \ref{inj} for a fixed $s\in\{0,\ldots, n\}$ we get
\begin{align*}
\tilde K_s(x,y) & = [A_{s,j} Q_{s+1}(x)]^\trans Q_s(y) - Q_s(x)^\trans [A_{s,j} Q_{s+1}(y)]
\\
& = [x_jQ_s(x)- B_{s,j}Q_s(x) - A_{s-1,j}^\trans Q_{s-1}(x)]^\trans Q_s(y)
\\
& \quad - Q_s(x)^\trans [y_jQ_s(y)- B_{s,j}Q_s(y) - A_{s-1,j}^\trans Q_{s-1}(y)]
\\
& \quad + r_s(x)^\trans Q_s(y) + Q_s(x)^\trans \tilde r_s(y),
\end{align*}
where $r_s$ and $\tilde r_s$ are column polynomials equal to $0$ modulo $V$. Remembering that $B_{k,j}$ are real and symmetric (see Remark \ref{inj}) we get 
\begin{align*}
\tilde K_s(x,y) & = (x_j-y_j) Q_s(x)^\trans Q_s(y)
\\
& \quad - [A_{s-1,j}^\trans Q_{s-1}(x)]^\trans Q_s(y) + Q_s(x)^\trans [A_{s-1,j}Q_{s-1}(y)]
\\
& \quad +r_s(x)^\trans Q_s(y) + Q_s(x)^\trans \tilde r_s(y)
\\
& = (x_j-y_j) Q_s(x)^\trans Q_s(y) + \tilde K_{s-1}(x,y)
+r_s(x)^\trans Q_s(y) + Q_s(x)^\trans \tilde r_s(y).
\end{align*} 
Technically, this procedure works only for $s\geqslant 1$ but it is easy to see that the term $\tilde K_{s-1}$ for $s=0$ is equal to zero, since $Q_{-1}=0$. Taking summation over $s=0,\ldots,n$ we complete the proof.
\end{proof}
The assertion of Theorem \ref{CDformula} can be written briefly as
\begin{align}\label{brief}
(x_j-y_j) K_n(x,y) \emvdwa [A_{n,j}Q_{n+1}(x)]^\trans Q_n(y) - Q_n(x)^\trans [A_{n,j}Q_{n+1}(y)].
\end{align}
Note that in case of $V=\nul$ the above equation coincides with the ordinary one. We say that $K_n$ defined by \eqref{cdk} {\em satisfies the $j$-th Christoffel-Darboux formula} if \eqref{brief} holds true with some scalar matrix $A_{n,j}$.
\begin{lem}\label{tensor}
Let $X$ and $Y$ be real or complex vector spaces, and $X_0$ and $Y_0$ be their linear subspaces, respectively. Then there is a linear isomorphism
$$
\varPhi: (X/X_0) \otimes (Y/Y_0) \to (X \otimes Y)/(X_0\otimes Y + X \otimes Y_0),
$$
such that $\varPhi((x+X_0)\otimes (y+Y_0)) = x\otimes y + X_0\otimes Y + X \otimes Y_0$ for all $x\in X$ and $y\in Y$.
\end{lem}
\proof
The authors are aware that this lemma is commonly regarded as a straightforward consequence of the Yoneda lemma, a well-established result in the category theory. However, we have not been able to find any direct reference in the mathematical literature (apart from fishy websites), hence we propose our short proof without recourse to any ``abstract nonsense.''

For convenience, we abbreviate $X_0\otimes Y + X \otimes Y_0$ to $\mathsf T(X_0,Y_0)$. We begin with showing that the mapping
$$
(X/X_0) \times (Y/Y_0) \ni (x+X_0,y+Y_0) \mapsto x\otimes y +  \mathsf T(X_0,Y_0) \in (X \otimes Y)/\mathsf T(X_0,Y_0)
$$
is well defined. Let $x_1+X_0 = x_2+ X_0$ and $y_1+Y_0 = y_2 + Y_0$. Then 
\begin{align*}
x_1\otimes y_1 = x_2\otimes y_2 + (x_1-x_2)\otimes y_2 + x_1 \otimes (y_1-y_2), 
\end{align*}
thus $x_1\otimes y_1$ and $x_2\otimes y_2$ are equal modulo $\mathsf T(X_0,Y_0)$ and the mapping is well defined. Since it is bilinear, the universal factorization property for tensor products yields the linear mapping
$$
\varPhi: (X/X_0) \otimes (Y/Y_0) \to (X \otimes Y)/\mathsf T(X_0,Y_0)
$$
resulting in $\varPhi ((x+X_0)\otimes (y+Y_0)) = x\otimes y +  \mathsf T(X_0,Y_0)$ for all $x\in X$ and $y\in Y$.

We now proceed to define a mapping which will promptly turn out to be the inverse of $\varPhi$. Consider the bilinear mapping
$$
X\times Y \ni (x,y) \mapsto (x+X_0)\otimes (y+Y_0) \in (X/X_0) \otimes (Y/Y_0).
$$
The bilinearity of this mapping and the universal factorization property lead to the linear mapping
$$
\varPsi : X\otimes Y \to (X/X_0) \otimes (Y/Y_0)
$$
such that $\varPsi (x\otimes y) = (x+X_0) \otimes (y+Y_0)$ for all $x\in X$ and $y\in Y$. It is clear that $\varPsi (x\otimes y) = 0$, whenever $x\in X_0$ or $y\in Y_0$. Since the elements of the form $x\otimes y$ with $x\in X_0$ or $y\in Y_0$ generate $\mathsf T(X_0,Y_0)$ it follows that $\mathsf T(X_0,Y_0) \subset \ker \varPsi$.  Passing to the quotient spaces the mapping $\varPsi$ induces another mapping
$$
\varPsi_1: (X\otimes Y)/\mathsf T(X_0,Y_0) \to (X/X_0) \otimes (Y/Y_0)
$$
with the property that $\varPsi_1 (x\otimes y + \mathsf T(X_0,Y_0)) = (x+X_0) \otimes (y+Y_0)$ for all $x\in X$ and $y\in Y$. It turns out that $\varPhi \circ \varPsi_1$ and $\varPsi_1 \circ \varPhi$ are identity mappings on $(X\otimes Y)/\mathsf T(X_0,Y_0)$ and $(X/X_0) \otimes (Y/Y_0)$, respectively, which can be verified directly on the generators.
\endproof
The following lemma shows that the ideal $V_2$ is well chosen regarding tensor products.
\begin{cor}\label{trans}
Let $V$ be a $*$-ideal in $\pcal_d$. There is a unique mapping
$$
\varPhi : (\pcal_d/V) \otimes (\pcal_d/V) \to \pcal_{2d}/V_2
$$
such that $\varPhi((p+V) \otimes (q+V)) = p\otimes q+V_2$. Moreover, $\varPhi$ is a linear isomorphism.
\end{cor}
\proof
This is Lemma \ref{tensor} translated to the case of a $*$-ideal treated as a linear subspace of $\pcal_d$. The identification of $\pcal_d \otimes \pcal_d$ with $\pcal_{2d}$ is well known.
\endproof
\begin{cor}\label{p-tensor}
Let $V\subset \pcal_d$ be a proper $*$-ideal and $\{Q_k\}_{k=0}^{\varkappa_V}$ be a rigid $V$-basis of $\pcal_d$. If $P_1$ and $P_2$ are column polynomials such that
\begin{equation}\label{cross}
P_1(x)^\trans Q_k(y) \emvdwa Q_j(x)^\trans P_2(y) 
\end{equation}
with some integer $k,j \geqslant 0$, then there exist a unique scalar matrix $E$ such that
$$
P_1 \emv E Q_j \quad \text{and} \quad P_2 \emv E^\trans Q_k.
$$
\end{cor}
\proof
It is well known that if $\{e_n\}_{n=0}^\varkappa$ with $\varkappa \in \nbb\cup \{\infty\}$  is a basis of a vector space $X$, then the system $\{e_m\otimes e_n\}_{m,n=0}^\varkappa$ forms a basis for $X\otimes X$. Hence, the system $\{Q_k\}_{k=0}^{\varkappa_V}$ induces the basis of $(\pcal_d/V) \otimes (\pcal_d/V)$ composed of all elements of the form $(p+V) \otimes (q+V)$, where $p$ appears in some column $Q_j$ and $q$ does in some column $Q_k$. By Corollary \ref{trans} the related basis of $\pcal_{2d}/V_2$ consists of all elements of the form $p\otimes q+V_2$ with the same way of choosing polynomials $p$ and $q$.

By our assumption $P_1$ can be expressed as $P_1 \emv \sum_{n=0}^N C_n Q_n$ with unique scalar matrices $C_n$. Similarly, $P_2 \emv \sum_{n=0}^N D_n Q_n$ with unique scalar matrices $D_n$ (adding some zero terms we may have the same $N$ for both $P_1$ and $P_2$). Substituting this in \eqref{cross} we get
\begin{align*}
\sum_{n=0}^N Q_n^\trans (x) C_n^\trans Q_k(y) \emvdwa \sum_{n=0}^N Q_j(x)^\trans D_n Q_n (y).
\end{align*}
Since the left hand side is a sum of polynomials of the form $p(x)q(y)$ with $q$ from the column $Q_k$, by linear independence modulo $V_2$ we infer that $D_n=0$ for all $n$ but $n=k$. Similarly, all $C_n=0$ except for $n=j$. It follows that $P_1= C_jQ_j$ and $P_2 = D_k Q_k$. Moreover, the equation \eqref{cross} takes the form 
$$
Q_j^\trans (x) C_j^\trans Q_k(y) \emvdwa Q_j^\trans (x) D_k Q_k(y),
$$
which by linear independence forces $C_j^\trans = D_k$. The proof is complete.
\endproof
As in one variable, we would like to emphesize that there is a deeper connection between the Christoffel-Darboux formula and the three term recurrence relation. In fact, the aforesaid formula implies the relation.
\begin{thm}
Let $V\subset\pcal_d$ be a proper $*$-ideal and $\{Q_k\}_{k=0}^{\varkappa_V}$ be a rigid $V$-basis of real polynomials with $Q_0=1$, and $j\in\{1,\ldots,d\}$. Let $N\in\nbb$, $N\leqslant \varkappa_V$. Assume that $K_s$ satisfies the $j$-th Christoffel-Darboux formula \eqref{brief} for all $n=0,\ldots,N$ with some scalar matrices $A_{n,j}$. Then 
$$
X_jQ_k \emv A_{k,j}Q_{k+1}+ B_{k,j}Q_k + A_{k-1,j}^\trans Q_{k-1}, \quad k=0,\ldots,N,
$$
with some scalar matrices $B_{k,j}$ $($with the convention $Q_{-1}=0$ and $A_{-1,j}=0)$.
\end{thm}
\proof
Let $\tilde K_n$ stand for the right hand side of \eqref{brief}, $n= 0,\ldots,N$. By our assumption
$$
\tilde K_n(x,y) = (x_j-y_j)K_n(x,y) + r_n(x)^\trans w_n(y) + \tilde w_n(x)^\trans \tilde r_n(y)
$$
with column polynomials $r_n$, $w_n$, $\tilde r_n$, $\tilde w_n$ such that $r_n$ and $\tilde r_n$ are equal to $0$ modulo $V$. For simplicity we write $R_n(x,y) = r_n(x)^\trans w_n(y) + \tilde w_n(x)^\trans \tilde r_n(y)$. Since 
\begin{equation*}
K_n(x,y) = Q_n(x)^\trans Q_n(y) + K_{n-1}(x,y),
\end{equation*}
we have
\begin{align*}
\tilde K_k(x,y) - R_k(x,y) = (x_j-y_j)Q_k(x)^\trans Q_k(y) + \tilde K_{k-1}(x,y) - R_{k-1}(x,y)
\end{align*}
for $k=0,\ldots,N$. This, written explicitly, reads as follows:
\begin{multline*}
[A_{k,j}Q_{k+1}(x)]^\trans Q_k(y) - Q_k(x)^\trans [A_{k,j}Q_{k+1}(y)] - R_k(x,y)
\\
= (x_j-y_j) Q_k(x)^\trans Q_k(y)+ [A_{k-1,j}Q_k(x)]^\trans Q_{k-1}(y)
\\
- Q_{k-1}(x)^\trans [A_{k-1,j}Q_k(y)] - R_{k-1}(x,y),
\end{multline*}
which, after rearranging terms, leads to
\begin{multline*}
[A_{k,j}Q_{k+1}(x) + A_{k-1,j}^\trans Q_{k-1}(x) - x_jQ_k(x)]^\trans Q_k(y) - R_k(x,y)
\\
= Q_k(x)^\trans [A_{k,j}Q_{k+1}(y) + A_{k-1,j}^\trans Q_{k-1}(y) - y_jQ_k(y)] - R_{k-1}(x,y).
\end{multline*}
A direct application of Corollary \ref{p-tensor} leads to a matrix $B_{k,j}$ such that
\begin{equation*}
A_{k,j}Q_{k+1} + A_{k-1,j}^\trans Q_{k-1} - x_jQ_k \emv - B_{k,j} Q_k
\end{equation*}
which gives the desired conclusion.
\endproof
\section{Examples}
One of examples could be simply rewritten from \cite{szaf1} where the author considered the tensor product of two families of orthogonal polynomials in one variable, more specifically the Krawtchouk polynomials and the Charlier ones. Digesting this example one may grasp a general background idea that allows to 
construct what may be called ``product polynomials'', i.e. the tensor product of orthogonal polynomials in a single variable. Noticing this we prefer to direct our attention to two variable polynomials which are not obtained this way.

We now focus on $\pcal_2/V$, where $V$ is the ideal of all polynomials vanishing on the unit circle $\tbb$, centered at $0$. In order to establish the three term recurrence relation we will employ the well known system $\{z^n\}_{n\in\zbb}$, which is orthonormal with respect to the normalized Lebesgue measure $m$ on $\tbb$. One can easily check that the system
\begin{equation}\label{syst}
\{z^0\} \cup \big\{ \tfrac 1{\sqrt 2} (z^n+\bar z^n)\colon n\geqslant 1 \big\} \cup \big\{ \tfrac 1{\sqrt 2\I}(z^n-\bar z^n) \colon n\geqslant 1 \big\}
\end{equation}
is also orthonormal with respect to $m$. What is more, every member of $\cbb[z,\bar z]$ is equal modulo $V$ to a linear combination of elements of the system. Indeed, since $z \bar z \emv 1$, for $k> l$ we may write
\begin{equation*}
z^k \bar z^l \emv z^{k-l} = \tfrac 12 (z^{k-l}+\bar z^{k-l}) + \tfrac 12 (z^{k-l}-\bar z^{k-l}),
\end{equation*}
so we have expressed $z^k\bar z^l$ linearly by means of the system \eqref{syst}. The same can be done for $k<l$, as $z^k\bar z^l \emv \bar z^{l-k}$. Since the case $k=l$ is trivial, this proves our claim.

These polynomials give rise to the following real polynomials:
\begin{align*}
& q^{(0)}_1 (x_1,x_2) =1,
\\&
\begin{gathered}
q^{(k)}_1 (x_1,x_2) = \tfrac 1{\sqrt 2} \big( (x_1+\I x_2)^k + (x_1-\I x_2)^k \big),
\\ 
q^{(k)}_2 (x_1,x_2) = \tfrac 1{\sqrt 2 \I} \big( (x_1+\I x_2)^k - (x_1-\I x_2)^k \big),
\end{gathered}
\quad k\geqslant 1.
\end{align*}
We arrange them in a column form:
\begin{equation}\label{column}
Q_0 = q^{(0)}_1, \quad Q_k = \begin{bmatrix} q^{(k)}_1 \\ q^{(k)}_2 \end{bmatrix}, \ k\geqslant 1
\end{equation}
Let $L_m:\pcal_2 \to \cbb$ denote the 
functional
\begin{equation*}
L_m(p) = \int_\mathbb T p(x_1,x_2) \D m(x_1,x_2), \quad p\in \pcal_2.
\end{equation*}
By the properties of the system \eqref{syst} and the identification of $\pcal_2$ with $\cbb[z,\bar z]$ via $z=x_1+\I x_2$ we infer that \eqref{column} is a rigid $V$-basis of $\pcal_2$ such that $L_m(Q_k Q_l^\trans) = 0$ if $k\neq l$, and $L_m(Q_k Q_k^\trans) = \big[ \begin{smallmatrix} 1&0\\0&1 \end{smallmatrix} \big]$ if $k\geqslant 1$. By Theorem \ref{gfav} the system $\{Q_k\}_{k=0}^\infty$ satisfies the three term recurrence relation
\begin{equation*}
X_jQ_k \emv A_{k,j}Q_{k+1}+ B_{k,j}Q_k + A_{k-1,j}^\trans Q_{k-1}, \quad k\geqslant 0,\ j=1,2.
\end{equation*}
Following the idea from Remark \ref{inj}, we may compute matrices $A_{k,j}$ and $B_{k,j}$ with the help of the formulas
\begin{align*}
A_{k,j} = L_m(X_jQ_kQ_{k+1}^\trans), \quad B_{k,j} = L_m(X_j Q_k Q_k^\trans), \quad k\geqslant 0,\ j=1,2.
\end{align*}
Leaving simple though slightly tedious computations to the reader one may arrive at
\begin{align*}
& A_{0,1}= \big[ \tfrac 1{\sqrt 2}\ 0 \big], \quad A_{0,2}= \big[ 0\ \tfrac 1{\sqrt 2} \big],
\\&
A_{k,1} = \begin{bmatrix} \frac 12 & 0 \\ 0 & \frac 12 \end{bmatrix}, \quad A_{k,2} = \begin{bmatrix} 0 & \frac 12 \\ -\frac 12 & 0 \end{bmatrix},\quad k\geqslant 1,
\end{align*}
and $B_{k,j} = 0$ for all $k\geqslant 0$, $j=1,2$. Theorem \ref{CDformula} is now ready to be applied for writing the Christoffel-Darboux formula
\begin{align*}
(x_j-y_j) \sum_{k=0}^n Q_k^\trans(x) Q_k(y)\emvdwa [A_{n,j} Q_{n+1}(x)]^\trans Q_n(y) - Q_n(x)^\trans [A_{n,j} Q_{n+1}(y)],
\end{align*}
where all the involved objects can be explicitly written.

Let us now turn to a more interesting example of the Bernoulli lemniscate $B$, i.e.\ the set of all points $(x_1,x_2) \in \rbb^2$ satisfying $(x_1^2+x_2^2)^2 = x_1^2-x_2^2$. In the case of polynomials in one complex variable this case was discussed in \cite{marc} (cf.\ \cite{marc2}). Our approach seems to be separated from that of \cite{marc}, not to mention no apparent relation between orthogonality of polynomials in complex and real variables (the precedent case of the circle has just been a lucky coincidence). A parametric description of $B$ is given by
$$
[-\tfrac \pi4, \tfrac \pi 4] \ni t \mapsto  \pm \gamma(t) \in \rbb^2, \quad \gamma(t) = 
\sqrt{\cos 2t} (\cos t, \sin t),
$$
which leads to the formula defining measure $\mfrak$ on $B$:
\begin{align*}
L_\mfrak (f) := \int_B f(x_1,x_2) \D \mfrak (x_1,x_2) & = \alpha_\mfrak \int_{-\frac \pi4} ^{\frac \pi 4} \big( f(\gamma(t)) + f(-\gamma(t)) \big) \|\gamma'(t)\| \D t
\\&
= \alpha_\mfrak \int_{-\frac \pi4} ^{\frac \pi 4} \big( f(\gamma(t)) + f(-\gamma(t)) \big) \frac{\D t}{\sqrt{\cos 2t}},
\end{align*}
where $f$ is a Borel function bounded on $B$ and $\alpha_\mfrak >0$ is chosen so that $\mfrak$ is normalized by $\int_B 1\D \mfrak =1$. It follows that $L_\mfrak(f) = 0$ whenever\label{abc}
\begin{itemize}
\item[(a)] $f(-x_1,-x_2)= -f(x_1,x_2)$ for all $(x_1,x_2)\in B$, or
\item[(b)] $f(-x_1,x_2)= -f(x_1,x_2)$ for all $(x_1,x_2)\in B$, or
\item[(c)] $f(x_1,-x_2)= -f(x_1,x_2)$ for all $(x_1,x_2)\in B$.
\end{itemize}
Indeed, (a) is a direct consequence of the formula for $\mfrak$, while (b) and (c) result from
$$
\int_{-\frac \pi4} ^{\frac \pi 4} f(\pm \gamma (t))  \frac{\D t}{\sqrt{\cos 2t}}  = \int_{-\frac \pi4} ^{\frac \pi 4} f(\pm \gamma (-t)) \frac{\D t}{\sqrt{\cos 2t}}.
$$
The ideal $V$ related to $B$ consists of all polynomials $p\in \pcal_2$ vanishing on $B$, or, equivalently, satisfying $L_\mfrak(|p|^2)=0$. In particular $X_1^2 - X_2^2 \emv (X_1^2 + X_2^2)^2$. This yields
\begin{equation}\label{kwadraty}
X_1^2 \emv \tfrac 12 (X_1^2+X_2^2 + (X_1^2 + X_2^2)^2) \quad \text{and} \quad X_2^2 \emv \tfrac 12 (X_1^2+X_2^2 - (X_1^2 + X_2^2)^2),
\end{equation}
which means that $X_1^{2j} X_2^{2k} \emv q(X_1^2+X_2^2)$ with some one-variable polynomial $q$ depending on nonnegative integers $j$ and $k$. As a consequence, if $p \in \pcal_2$, then there exist $q_j \in \pcal_1$, $j=1,2,3,4$, such that 
\begin{equation*}
p \emv q_1(X_1^2+X_2^2) + X_1 q_2(X_1^2+X_2^2) + X_2 q_3(X_1^2+X_2^2) + X_1X_2 q_4(X_1^2+X_2^2).
\end{equation*}
The crucial remark to make here is that every $p\in \pcal_2$ is equal modulo $V$ to a linear combination of polynomials  $p^{j,k}_l := X_1^j X_2^k (X_1^2 + X_2^2)^l$ with integers $j,k=0,1$ and $l\geqslant 0$. It follows that if $j_1,j_2,k_1,k_2 =0,1$, $(j_1,k_1) \neq (j_2,k_2)$ and $q_1,q_2\in\pcal_1$, then 
\begin{equation}\label{x1j1}
X_1^{j_1}X_2^{k_1} q_1(X_1^2+X_2^2) \cdot X_1^{j_2}X_2^{k_2} q_2(X_1^2+X_2^2) \emv 
X_1^r X_2^s q(X_1^2+X_2^2),
\end{equation}
with some $q \in \pcal_1$ and $r,s =0,1$ such that $r+s>0$. Since the right hand side of \eqref{x1j1} is a function satisfying one of the properties (a), (b) and (c) listed on p.\ \pageref{abc}, we infer that
$$
L_\mfrak \big( X_1^{j_1}X_2^{k_1} q_1(X_1^2+X_2^2) \cdot X_1^{j_2}X_2^{k_2} q_2(X_1^2+X_2^2) \big) = 0
$$
with all $j_1,j_2,k_1,k_2,q_1,q_2$ already chosen. In other words, polynomials 
$$
X_1^{j_1}X_2^{k_1} q_1(X_1^2+X_2^2) \quad \text{and} \quad X_1^{j_2}X_2^{k_2} q_2(X_1^2+X_2^2)
$$
are {\em $L_\mfrak$-orthogonal}. In particular, polynomials $p^{j_1,k_1}_{l_1}$ and $p^{j_2,k_2}_{l_2}$ are $L_\mfrak$-orthogonal, provided $j_1,j_2,k_1,k_2$ are as above and $l_1,l_2 \geqslant 0$.

Fix any $j,k=0,1$ and consider the linear functional
$$
L^{j,k} : \pcal_1 \to \cbb, \quad L^{j,k}(q) = L_\mfrak(X_1^{2j} X_2^{2k} q(X_1^2+X_2^2)), \ q\in \pcal_1.
$$
It is evident that $L^{j,k}(|q|^2) >0$ for any $q \neq 0$, thus there exists a sequence of real orthogonal polynomials $\{q_l^{j,k}\}_{l=0}^\infty$ such that $\deg q_l^{j,k} = l$, $l \geqslant 0$, and $$
L^{j,k}(q_{l_1}^{j,k} q_{l_2}^{j,k}) = \delta_{l_1l_2}.
$$
As usual for orthogonal polynomials in one variable the sequence is determined uniquely up to the unimodular factor of each $q_l^{j,k}$. It is a matter of routine to verify that the family of all polynomials 
\begin{equation}\label{L-ort}
W^{j,k}_l := X_1^j X_2^k q_l^{j,k}(X_1^2 + X_2^2), \quad j,k=0,1, \ l\geqslant 0,
\end{equation}
is $L_\mfrak$-orthogonal, and even $L_\mfrak$-orthonormal. We now arrange this family in the column form:
$$
Q_0 = 1, \ 
Q_1 = \begin{bmatrix} W^{1,0}_0 \\[.5ex] W^{0,1}_0\end{bmatrix}\hspace{-.7ex}, \ 
Q_2 = \begin{bmatrix} W^{0,0}_1 \\[.5ex] W^{0,0}_2 \\[.5ex] W^{1,1}_0
\end{bmatrix}\hspace{-.7ex}, \
Q_{2s-1} = \begin{bmatrix} W^{1,0}_{2s-3} \\[.5ex] W^{1,0}_{2s-2} \\[.5ex] W^{0,1}_{2s-3} \\[.5ex] W^{0,1}_{2s-2}\end{bmatrix}\hspace{-.7ex},\
Q_{2s} = \begin{bmatrix} W^{0,0}_{2s-1} \\[.5ex] W^{0,0}_{2s} \\[.5ex] W^{1,1}_{2s-3} \\[.5ex] W^{1,1}_{2s-2}\end{bmatrix}\hspace{-.7ex},\ s\geqslant 2.
$$
At first glance this does not look like a good idea leading to a rigid $V$-basis of $\pcal_2$, because some coordinate polynomials in $Q_n$ are not of degree $n$ if $n\geqslant 2$. It turns out that for all $n\geqslant 0$ every coordinate polynomial of $Q_n$ is equal modulo $V$ to a polynomial of degree exactly $n$. This follows from the equation\footnote{\hspace{.3ex}see the notation introduced in the beginning of Section \ref{CDreal}}
\begin{equation}\label{piv}
\Pi_V (\pcal^k_2) = \lin \Pi_V \Big( \bigcup_{n=0}^k Q_n \Big) \quad\text{for every } k\geqslant 0,
\end{equation}
which we are now going to prove by induction. The instances of $k=0$ and $k=1$ are trivial, while for $k=2$ we have $\deg W^{0,0}_1 = \deg W^{1,1}_0 = 2$ and $W^{0,0}_2$ is equal modulo $V$ to a polynomial of degree $2$, since $(X_1^2+X_2^2)^2 \emv X_1^2-X_2^2$. Let us now assume that $k\geqslant 3$. To prove the inclusion ,,$\subset$'' fix a monomial $X_1^\alpha X_2^\beta$ such that $\alpha + \beta \leqslant k$. By induction hypothesis we may focus only on the case when $\alpha + \beta =k$. Let $[m]_2 = 0$ if integer $m$ is even, and $[m]_2=1$ if $m$ is odd. Since 
$$
X_1^\alpha X_2^\beta = X_1^{[\alpha]_2} X_2^{[\beta]_2} X_1^{\alpha - [\alpha]_2}X_2^{\beta - [\beta]_2},
$$
in virtue of \eqref{kwadraty} we get
\begin{equation*}
X_1^\alpha X_2^\beta \emv X_1^{[\alpha]_2} X_2^{[\beta]_2} p_{\alpha,\beta} (X_1^2+X_2^2),
\end{equation*}
with a polynomial $p_{\alpha,\beta}$ such that $\deg p_{\alpha,\beta} = \alpha+\beta - [\alpha]_2 - [\beta]_2$. Under the notation $\mu= [\alpha]_2$ and $\nu = [\beta]_2$ expanding $p_{\alpha,\beta}$ in the basis $\{q_l^{\mu,\nu}\}_{l=0}^\infty$ we obtain
$$
X_1^\alpha X_2^\beta \emv X_1^\mu X_2^\nu \sum_{l=0}^{k-\mu-\nu} a_l q_l^{\mu,\nu}(X_1^2+X_2^2)
= \sum_{l=0}^{k-\mu-\nu} a_l W_l^{\mu,\nu}
$$
with some coefficients $a_l\in \cbb$. Considering all four possible cases of $\mu,\nu =0,1$ we now show that every $W_l^{\mu,\nu}$ in the above sum is a coordinate in one of $Q_n$, $n=0,\ldots,k$. Indeed, if $k=2s-1$, then $(\mu,\nu)$ is equal to $(1,0)$ or $(0,1)$, so the resulting polynomials $W^{\mu,\nu}_l$ with $l=0,\ldots, 2s-2$ are coordinate polynomials in one of $Q_1,Q_3,\ldots,Q_{2s-1}$. In turn, if $k=2s$, then $(\mu,\nu)$ is equal to $(0,0)$ or $(1,1)$. In the case when $(\mu,\nu)=(0,0)$ we obtain $W^{0,0}_l$ with $l=0,\ldots,2s$ which appear in the columns $Q_0, Q_2, \ldots,Q_{2s}$. The remaining case $(\mu,\nu)=(1,1)$ can be done in a similar way. This proves the desired inclusion. To verify the reverse inclusion fix $\mu,\nu=0,1$ and take any coordinate polynomial $W_l^{\mu,\nu}$ in $Q_k$. By induction hypothesis it suffices to show that $W_l^{\mu,\nu} +V \in \Pi_V(\pcal_2^k)$. The polynomial $q_l^{\mu,\nu}$ can be written as $q_l^{\mu,\nu} = \sum_{j=0}^l b_j X^j$ with some coefficients $b_l\in \cbb$ and therefore
\begin{equation}\label{jedno}
W_l^{\mu,\nu} = X_1^\mu X_2^\nu q_l^{\mu,\nu} (X_1^2+X_2^2) = \sum_{j=0}^l b_j X_1^\mu X_2^\nu (X_1^2+X_2^2)^j.
\end{equation}
Since $W_l^{\mu,\nu}$ is a coordinate of $Q_k$, we see that $\mu+\nu +l \leqslant k$. By \eqref{jedno} it suffices to show that
$$
X_1^\mu X_2^\nu (X_1^2+X_2^2)^j + V \in \Pi_V(\pcal_2^k), \quad \mu+\nu+j \leqslant k.
$$
As $(X_1^2+X_2^2)^2 \emv X_1^2-X_2^2$, for an even number $j$ we get
$$
X_1^\mu X_2^\nu (X_1^2+X_2^2)^j \emv X_1^\mu X_2^\nu (X_1^2-X_2^2)^{j/2},
$$
and the latter polynomial belongs to $\pcal_2^k$. If in turn $j$ is odd we have
\begin{align}\label{munu}
\begin{aligned}
X_1^\mu X_2^\nu (X_1^2+X_2^2)^j & = X_1^\mu X_2^\nu (X_1^2+X_2^2) (X_1^2+X_2^2)^{j-1}
\\&
\emv X_1^\mu X_2^\nu (X_1^2+X_2^2) (X_1^2-X_2^2)^{(j-1)/2},
\end{aligned}
\end{align}
where we are led to polynomial of degree $\mu+\nu + j + 1$, seemingly exceeding $k$. However, if $j<l$, then the last polynomial in \eqref{munu} is of degree less than or equal to $k$. In the remaining case when $j=l$ is odd one may notice that whenever $W^{\mu,\nu}_l$ appears in $Q_k$, then $\mu+\nu+l<k$, and the resulting polynomial in \eqref{munu} is again a member of $\pcal_2^k$. We have thus proved \eqref{piv}.

Since all $W_l^{j,k}$ are $L_\mfrak$-orthonormal,  we deduce that they are linearly independent modulo $V$, so by \eqref{piv} the elements $W_l^{j,k}+V$ form a basis of $\pcal_2/V$. It follows that the system $\{Q_n\}_{n=0}^\infty$ meets all the requirements to be a rigid basis except for the condition on the degree of coordinate polynomials of $Q_n$. However, this can be easily overcome by the above discussion where we have shown how the coordinate polynomials can be replaced by polynomials of the proper degree preserving equality modulo $V$. Since our goal is to write the three term recurrence relation as in Theorem \ref{gfav}, which is equation modulo $V$, i.e.
$$
X_jQ_k \emv A_{k,j}Q_{k+1}+ B_{k,j}Q_k + A_{k-1,j}^\trans Q_{k-1},\quad k\in\nbb,\ j=1,2,
$$
we do not really have to bother about the degree of polynomials in $Q_n$. As we already know,
$$
A_{k,j} = L_\mfrak(X_jQ_kQ_{k+1}^\trans) \quad \text{and} \quad B_{k,j} = L_\mfrak(X_jQ_kQ_k^\trans).
$$
A typical entry in $A_{k,j}$ or $B_{k,j}$ takes the form $L_\mfrak(X_j W_{l_1}^{\mu_1,\nu_1} W_{l_2}^{\mu_2,\nu_2})$ which is equal to zero but the following two cases:
\begin{enumerate}
\item[(i)] $j=1$, $\mu_1+\mu_2=1$, $\nu_1=\nu_2$, or
\item[(ii)] $j=2$, $\mu_1=\mu_2$, $\nu_1 + \nu_2 =1$,
\end{enumerate}
which can be easily justified with the help of properties (a), (b) and (c) on p.\ \pageref{abc}. Since all the polynomials $W^{\mu,\nu}_l$ in the column $Q_k$ satisfy either $\mu+\nu=1$ (if $k$ is odd) or $\mu+\nu \in \{0,2\}$ (if $k$ is even), none of the entries of the matrix $X_jQ_kQ_k^\trans$ satisfies (i) or (ii), hence $B_{k,j} =0$ for all admissible $k$ and $j$. We now consider $A_{k,j}$. Fix $s\geqslant 2$ and introduce the auxiliary column polynomial $R_l^{\mu,\nu} = [q_{l-1}^{\mu,\nu} \ q_l^{\mu,\nu} ]^\trans$, $l \geqslant 1$, $\mu,\nu=0,1$. 
Thus we may write
\begin{align*}
Q_{2s-1} = \begin{bmatrix} X_1R^{1,0}_{2s-2}( X_1^2 + X_2^2) \\[1ex] X_2 R^{0,1}_{2s-2} ( X_1^2 + X_2^2)\end{bmatrix}
\quad \text{and} \quad
Q_{2s} = \begin{bmatrix} R^{0,0}_{2s}( X_1^2 + X_2^2)
\\[1ex] X_1X_2 R^{1,1}_{2s-2} ( X_1^2 + X_2^2)\end{bmatrix}, \ s\geqslant 2,
\end{align*}
where $R^{\mu,\nu}_l ( X_1^2 + X_2^2)$ means substituting $X_1^2+X_2^2$ in the both coordinate polynomials of $R^{\mu,\nu}_l$. In virtue of (i) and (ii) we get
\begin{align*}
A_{2s,1} & = \begin{bmatrix}
L^{1,0} (R^{0,0}_{2s} (R^{1,0}_{2s})^\trans) & 0\\
0 & L^{1,1}(R^{1,1}_{2s-2} (R^{0,1}_{2s})^\trans)
\end{bmatrix},
\\
A_{2s-1,1} & = \begin{bmatrix}
L^{1,0} (R^{1,0}_{2s-2} (R^{0,0}_{2s})^\trans) & 0\\
0 & L^{1,1}(R^{0,1}_{2s-2} (R^{1,1}_{2s-2})^\trans)
\end{bmatrix},
\\
A_{2s,2} & = \begin{bmatrix}
0 & L^{0,1} (R^{0,0}_{2s} (R^{0,1}_{2s})^\trans)\\
L^{1,1}(R^{1,1}_{2s-2} (R^{1,0}_{2s})^\trans) & 0
\end{bmatrix},
\\
A_{2s-1,2} & = \begin{bmatrix}
0 & L^{1,1} (R^{1,0}_{2s-2} (R^{1,1}_{2s-2})^\trans)\\
L^{0,1}(R^{0,1}_{2s-2} (R^{0,0}_{2s})^\trans) & 0
\end{bmatrix}.
\end{align*}
We encourage the reader to derive similar formulas for $A_{k,j}$ with $k=0,1$ and $j=0,1$, the cases not covered above. If we now employed  a one-variable orthonormalization procedure, we could consecutively compute all the matrices $A_{k,j}$.

It is worth noting that one may derive the integral representation for all the functionals $L^{j,k}$, $j,k=0,1$, that is
\begin{align*}
L^{j,k}(q) = \alpha_\mfrak \int_0^1 q(x) \frac{2^{1-j-k} x^{j+k} (1+x)^j (1-x)^k}{\sqrt{x(1-x^2)}} \D x, \quad q\in \pcal_1.
\end{align*}
The weight function for $L^{0,0}$ can be related to the weight mentioned in \cite[p. 37]{szego} (formula (2.9.1)), in order to see this it is enough to perform change of variable $t=1-x$ under the integral. In turn, the orthogonal polynomials associated with $L^{j,k}$ with $j,k=0,1$, $j+k>0$, can be derived from those of $L^{0,0}$ by Theorem 2.5 in \cite{szego}. 

The case considered in \cite{marc} concerns complex orthogonality of polynomials in a single complex variable with respect to an arbitrary admissible functional $L$ (let the term ``admissible'' remain mysterious for the time being). This allows us to introduce this case in a little bit sketchy way. We will show that the system 
\begin{equation}\label{sigma}
\Sigma_0 = 1,
\quad \Sigma_1 = \begin{bmatrix}X_1 \\ X_2\end{bmatrix},
\quad \Sigma_2 = \begin{bmatrix}X_1^2 \\ X_1X_2 \\ X_2^2 \end{bmatrix},
\quad \Sigma_n = \begin{bmatrix}X_1^n \\ X_1^{n-1} X_2 \\ X_1^{n-2} X_2^2 \\ X_1^{n-3} X_2^3 \end{bmatrix}, \ n\geqslant 3,
\end{equation}
forms a rigid $V$-basis of $\pcal_2$. By the discussion in \cite[Section 4]{css} the system $\{\Sigma_n\}_{n=0}^\infty$ is a proper candidate for a rigid $V$-basis of $\pcal_2$, because its structure is the same as that of the rigid $V$-basis $\{Q_n\}_{n=0}^\infty$ constructed above (i.e. the lengths of consecutive columns in $\{\Sigma_n\}_{n=0}^\infty$ are equal to  $d_V(n)$). Hence, it suffices to show that for every $j,k\geqslant 0$ we have
\begin{equation}\label{mono}
\Pi_V(X_1^j X_2^k) \in \lin \Pi_V \Big( \bigcup_{l=0}^{j+k} \Sigma_l \Big).
\end{equation}
If $j+k\leqslant 3$, then \eqref{mono} is satisfied in an apparent way, because then $X_1^j X_2^k$ is a member of $\Sigma_{k+j}$. Assume that $j+k \geqslant 4$ and $k \geqslant 4$, which covers all the monomials of degree $j+k$ outside of $\Sigma_{j+k}$. Since $X_1^2-X_2^2 \emv (X_1^2+X_2^2)^2$ we infer that
$$
X_2^4 \emv X_1^2 - X_2^2 -X_1^4 - 2X_1^2 X_2^2.
$$
This implies that
$$
X_1^j X_2^k  = X_1^j X_2^{k-4} X_2^4 \emv 
X_1^{j+2} X_2^{k-4}
-X_1^j X_2^{k-2}
-X_1^{j+4} X_2^{k-4}
-2 X_1^{j+2} X_2^{k-2}.
$$
This way we have reduced the powers of the variable $X_2$ by $2$. If $k-2 \geqslant 4$ we may repeat this procedure for all the monomials with powers of $X_2$ greater than or equal to $4$, eventually obtaining the linear combination of monomials with power of $X_2$ less than or equal to $3$. Thus we have shown \eqref{mono}. 

We now want to assume that $L$ is induced by a Borel measure $\mu$ whose support is a Zariski dense subset of the lemniscate $B$, i.e.\ for every $p\in \pcal_2$ if $p$ vanishes on the support of $\mu$, then it vanishes on the whole $B$ (for the Zariski topology and related notions see the discussion in \cite[Section 9]{css}). There is no loss of generality since by \cite[Theorem 43]{css} (see also \cite[Theorem 1]{schm}) for any positive definite functional $L:\pcal_d\to\cbb$ satisfying $L(q\bar q)=0$ with some $q\in\pcal_d$ such that the zero set $\mathcal Z_q := q^{-1}(\{0\})$ is compact\footnote{\ we allow the empty set here} in $\rbb^d$ there exists a Borel measure $\mu$ supported in $\mathcal Z_q$ such that
\begin{equation}\label{repr}
L(p) = \int_{\mathcal Z_q} p(x)\D\mu(x), \quad p\in \pcal_d. 
\end{equation}
This means in particular that if $L:\pcal_2 \to \cbb$ is positive definite and $L(p)=0$ for all $p\in\pcal_2$ such that $p|_B\equiv 0$, then $L$ must have integral representation \eqref{repr} with $\mathcal Z_p=B$.

It is worth noting that a subset of $B$ is Zariski dense in $B$ if and only if it is infinite. Indeed, it has to be infinite because every finite set is closed in the Zariski topology. In turn, assume that $B_1$ is an infinite subset of $B$, and $p\in\pcal_2$ vanishes on $B_1$. Employing the parametrizaton
$$
\varphi: \rbb \ni t \mapsto \Big( \frac{\cos t}{1+\sin^2 t}, \frac{\sin t \cos t}{1+\sin^2 t} \Big) \in B
$$
we see that $p \circ \varphi$ is real-analytic and the set of its zeros has an acumulation point, because there are infinitely many $t \in [0,2\pi]$ such that $\varphi(t)\in B_1$. By the identity principle for real-analytic functions (see \cite[Corollary 1.2.7]{krantz}) we infer that $p\circ \varphi$ is equal to zero everywhere, hence $p$ vanishes on $B$. 

Having fixed $L$ induced by a measure whose support is a Zariski dense subset of $B$ we may now perform an orthonormalization procedure applied to the system \eqref{sigma}. As usual, this can be done recursively, but it would be convenient to notice that once we have obtained columns $Q_0$, \ldots, $Q_{k-1}$ of orthonormal column polynomials (with lengths according to the structure of the system \eqref{sigma}), the column polynomial
$$
\widehat Q_k = \Sigma_k - \sum_{j=0}^{k-1} L(\Sigma_k Q_j^\trans)Q_j
$$
is orthogonal to all $Q_0$, \ldots, $Q_{k-1}$, hence it is orthogonal to $\pcal_2^{k-1}$. Thus  we may focus on orthonormalizing coordinate polynomials of $\widehat Q_k$. One of the possible ways of achieving this is finding a real matrix $U$ such that $U L(\widehat Q_k \widehat Q_k^\trans) U^\trans$ is the identity matrix, which is possible as $L(\widehat Q_k \widehat Q_k^\trans)$ is a nonsingular positive definite matrix (see \cite[Proposition 13]{css}). Then the formula $Q_n = U \widehat Q_n$ gives the desired column polynomial.

\medskip
{\bf Acknowledgment.} The first named author acknowledges the financial support of the Priority Research Area SciMat under the program Excellence Initiative Research University at the Jagiellonian University in Krak\'ow, Poland, decision no. U1U/P05/NO/03.55.

\bibliographystyle{amsalpha}

\end{document}